\documentclass{article}
\usepackage{xcolor}
\usepackage{amsmath}
\usepackage{amssymb}
\usepackage{graphicx}
\usepackage{comment}
\usepackage{setspace}
\usepackage{soul}
\usepackage{pgf}
\usepackage{tikz}
\usepackage{hyperref,mathrsfs}
\usepackage{arydshln}
\usetikzlibrary{shapes}
\newtheorem{theorem}{Theorem}[section]
\newtheorem{lemma}[theorem]{Lemma}
\newtheorem{proposition}[theorem]{Proposition}
\newtheorem{corollary}[theorem]{Corollary}
\newtheorem{remark}[theorem]{Remark}
\newtheorem{example}[theorem]{Example}
\newenvironment{proof}

\newcommand{\ct}[1]{#1^{\mbox{\tiny\sf T}}}            
\newcommand{\N}{{\mathbb N}}                           
\newcommand{\R}{{\mathbb R}}                           
\newcommand{\BC}{{\mathbb C}}                          
\newcommand{\BH}{{\mathbb H}}                          
\newcommand{\BD}{{\mathbb D}}                          
\newcommand{\BL}{{\mathbb L}}                          
\newcommand{\cA}{\mathscr{A}}                          
\newcommand{\cC}{\mathscr{C}}                          
\newcommand{\cD}{\mathscr{D}}                          
\newcommand{\cH}{\mathscr{H}}                          
\newcommand{\cM}{\mathscr{M}}                          
\newcommand{\cS}{\mathscr{S}}                          
\def\wek#1{\textsf{$\textbf{#1}$}}                     
\def\supwek#1{^{\mbox{\footnotesize $\wek #1$}}}       
\def\son{1}                                         
\def\sp{\wek p}
\def\sq{\wek q}

\def\dd{{\wek d}}
\def\sd{{\wek d}}

\def\sce{\wek c}
\def\si{{\wek 1}}
\def\se{\wek e}

\def\sd{\wek d}
\def\sx{\wek x}
\def\su{\wek u}

\def\sx{\wek x}
\def\sz{\wek 0}
\def\sxr{[\son+\sx\sr]}

\def\Ri#1{[\son+#1]}
\def\sy{\wek y}
\def\sr{\wek r}
\def\sR{\wek R}
\def\sD{\wek D}
\def\sN{\wek N}
\def\ort{^{\mbox{\tiny $\perp$}}}
\def\ct#1{#1^{\mbox{\tiny\sf T}}}

\def\ddmat#1#2#3#4{                                       
\begin{pmatrix}
#1&#2\\
#3&#4
\end{pmatrix} }
\def\dlog#1#2{\left(\genfrac{}{}{0pt}{}{#1}{#2}\right]}   
\def\bm#1{\mbox{\boldmath $#1$}}
\newcommand{\la}{\langle}
\newcommand{\ra}{\rangle}
\newcommand{\Proof}{\begin{proof}}
\newcommand{\QED}{\end{proof}}
\newcommand{\be}{\begin{equation}}
\newcommand{\ee}{\end{equation}}
\def\set#1{\left\{\,#1\,\right\} }                   
\newcommand{\df}{\,\stackrel{\mbox{\footnotesize\sf def}}{=}\,}
\newcommand{\imp}{\,\Rightarrow\,}
\newcommand{\tss}[1]{\textsuperscript{#1}}
\def\ca#1{{\cal #1}}
\def\MINT#1{\raisebox{-8pt}[10pt]{\begin{tabular}{c}
$\Int\!\!\cdots\!\!\Int$\\$^{#1}$\\\end{tabular}}}

\newcommand{\vspandex}{\rule[-18pt]{0pt}{36pt}}      
\newcommand{\vspandexsmall}{\rule[-11pt]{0pt}{22pt}} 
\newcommand{\vv}{

}                                                    

\newcommand{\E}{\mbox{\sf E}\,}                           
\def\ce#1#2{\,\E\left[\,#1\,|\,#2\,\right]\,}             
\newcommand{\Min}{\displaystyle\min}
\newcommand{\Prod}{\displaystyle\prod}
\newcommand{\Int}{\displaystyle\int}
\newcommand{\Sum}{\displaystyle\sum }
\newcommand{\Frac}{\displaystyle\frac}
\newcommand{\refrm}[1]{{\rm(\ref{#1})}}              
\newcommand{\tr}{{\rm tr}\,}                         
\begin{document}
\begin{center}
{\bf\large Algebraic structures spanned by differential-like operators 

on toy Fock spaces.}
\vspace{20pt}

{Jerzy Szulga}
\vspace {10pt}

{Department of Mathematics and Statistics,

Auburn University, Auburn, AL 36849, USA

Email: szulgje@auburn.edu}
\vspace{20pt}

\end{center}

\noindent {\bf Abstract. } {We study multiplicative systems of linear mappings acting on the toy Fock space, a.k.a.\  Rademacher chaos or  Walsh-Fourier series,  related to the creation, annihilation, and conservation operators in quantum probability. Like differential operators they entail analogs of the Leibnitz Formula and the Chain Rule, derived with the help of Riesz products (or discrete coherent vectors). Two symmetries among these operators entail groups and algebras of varying  signatures and mixed commutativity, generalizing Pauli spin matrices and quaternions. In particular, anticommutative Rademacher systems are constructed.}

\vspace{10pt}

\noindent {\bf Keywords:} {Toy Fock space, Riesz products,  creation-annihilation-conservation, discrete derivatives, rigged Hilbert spaces}

\vspace{10pt}
\noindent {\bf AMS 2010 MSC:} {Primary 47L30; Secondary 46L53, 60H07, 42C10}
\vspace{20pt}

\section{Introduction}
Walsh functions on the interval $[0,1]$, i.e., products of Rademacher functions \cite{Rad}, form an orthogonal basis of $\BH=L^2[0,1]$, and yield a substantial theory of Walsh-Fourier series \cite{SWS} (see also the recent survey \cite{Wei}), rich in applications and connections to other areas of mathematics. From another point of view Walsh functions provide an example of a commutative multiplicative group, which in turn entails the algebra  of ``discrete random chaos''.  The purpose of this note is to examine linear operators acting on this algebra of which many are analogs of differential operators. They begin with creation, annihilation, and conservation operators that were introduced and investigated by P-A. Meyer \cite{Mey} in the context of so called ``toy Fock spaces''. In particular, we derive suitable versions of the Leibnitz Formula and the Chain Rule.  \vv

These operators correspond to analogous operators (cf.\  K.R.\ Parthasarathy \cite{Par}) acting on a Fock space, which is essentially Gaussian chaos. The discrete context entails their simplified but not simplistic versions. Discrete approximation of Fock space naturally attracted interest of many authors (cf., e.g., \cite{Att,AttNech,Bel,Pau}). Various aspects of discrete stochastic calculus were discussed in \cite{Pri,NouPecRei} (both with a substantial bibliography).
Similar topics were touched by this author \cite[in Chaps. 8 \& 12]{Szu}, and some of them now are being expanded. \vv

In Section 2 we will briefly explain the notion of the toy Fock space and then modify Meyer's original set-theoretic notation to facilitate the algebraic and logical formalism of the forthcoming actions. We adopt the well known convention of generalized powers, e.g., used in \cite{Rud} to handle partial differential operators. In particular, the binary 0-1 framework makes most of results routinely programmable. \vv

Meyer's chaos operators can be interpreted in the context of discrete differentiation and integration, leading to analogs of the Leibnitz formula and Chain Rule. These and related topics are discussed in Section 3.\vv

In Section 4, from the plethora of chaos operators we select just two, related to the basic symmetries of the square, that can be seen as the first and the third Pauli spin matrices \cite[Sect.\ 2.1.3.]{Mey}. Denoted in this note by $R$ and $S$, these operators entail a seemingly simple yet rich multiplicative structure including generalizations of the notions of the spin matrices and quaternions. In particular, we give a matrix representation of the derived ``signed group'', as we call it, as well as an alternative representation that we name the ``double logic''. The latter structure is based on double 0-1 sequences that can be augmented freely, either at will, or according to prescribed patterns, or even at random. Again, due to the 0-1 framework, the augmentation follows a straightforward algorithm that controls the expanding eigen-system. Their matrix representations, or of products of the basic symmetries $R$ and $S$, naturally generate so called Clifford algebras, thus connecting them to quantum physics. All underlying operators, including Meyer's quantum operators, belong to these algebras.
\vv

In Section 5 we show that this structure  contains the anticommutative Rademacher systems which are constructed with the help of so called rigged Hilbert spaces.
\section{Framework}
\subsection{Toy Fock spaces}
Consider a connected subset $T\subset\R$. Let $(X_t,t\in T)$ be a standard Brownian motion, $\sigma(X)$ denote the $\sigma$-field generated by the increments $X_t-X_s$. Let $\BH=L^2(\sigma(X))$ be the Hilbert space of square-integrable $\sigma(X)$-measurable functions, called {\em functionals of Brownian Motion}. They admit a ``chaos representation'' \cite{Ito}:
 \[
 \sum_{k=0}^\infty X^k(f_k),\quad \mbox{where}\quad X^k(f_k)=\MINT {T^k} f(t_1,\dots,t_k)\,X(dt_1)\cdots X(dt_k),
 \]
with square integrable $k$-variate symmetric functions $f_k$, vanishing on diagonals. By convention, $X^0f_0\in\R$. We thus observe the direct sum $\BH=\bigoplus_k \BH_k$ of Hilbert spaces $\BH_k$. Its every element can be identified with a sequence $\wek f=(f_k)$. Under the norm
\[
 \|\wek f\|^2=\sum_n \frac{1}{k!}\|f_k\|_k^2.
 \]
the Hilbert space is called the {\em Fock space}. The ``$k$-homogeneous chaos'', or simply $\BH_k$, may be interpreted as a system of elementary particles, chaotically moving and interacting with each other. Further, the system may lose some particles or be enlarged by new ones. The dynamics and evolution of the chaotic systems can be modeled by a plethora of various operators acting within and between the chaos spaces, of which the Ornstein-Uhlenbeck semigroup of operators is a fundamental example. It suffices to define it on each $\BH_k$, at the same time recording its eigenstructure:
\[
T_t X^k f_k=e^{-kt} X^k f_k.
\]
Further, considering simple functions, the discrete forms of Gaussian chaos, spanned by independent $N(0,1)$ random variables,
\[
 \sum_{k=0}^\infty \,\,\sum_{n_1<...<n_k}\,a^{(k)}_{n_1,...,n_k}\,\gamma_{n_1}\cdots \gamma_{n_k}
 \]
yield a dense subset of $\BH$. Here we consider a finite or countable partition $(t_n)$ of $T$ and denote by $\gamma_n$ the $n$\tss{th} normalized increment:
\[
\gamma_n=\frac{X_{t_n}-X_{t_{n-1}}}   {\sqrt{|t_n-t_{n-1}|}}.
\]
Thus, the Ornstein-Uhlenbeck operator $T_t$ is sufficiently defined on the products
\be\label{OEG}
T_t\gamma_{n_1}\cdots \gamma_{n_k}=e^{-kt}\gamma_{n_1}\cdots \gamma_{n_k}.
\ee
While $T_t$ are bounded for $t\ge 0$, they are unbounded operators for $t<0$, so is the semigroup's generator
\be\label{Mall}
 A=\lim_{t\searrow 0} \frac{1}{t}(T_t-I),\qquad
 A \gamma_{n_1}\cdots \gamma_{n_k}=-k \,\gamma_{n_1}\cdots \gamma_{n_k}.
 \ee
Metaphorically speaking, the operators ``recognize'' the signs and are ``insensitive'' to scales. Indeed, speaking rigorously, signs are sufficient for their definition. In other words, conditioning upon signs $r_n$ of Gaussian random variables $\gamma_n$, which are independent of their moduli $|\gamma_n|$, we arrive at operators acting on a the discrete random chaos, spanned by the Rademacher functions. Thus, $L^2[0,1]$ with its orthonormal Walsh basis becomes a ``toy Fock space''. Of course, the Gaussian origins are not necessary, yet they build a strong expectation of potential applications  to the actual Fock spaces,  as shown by the references cited in Introduction.

\subsection{Products and sums}\label{notat}
Sequences are marked by the boldface upright font. For two special sequences we may use the double notation $\wek 1=\son=(1,1,1,...)$ and $\sz=0=(0,0,0,\dots)$ (the distinction between the number and the sequence will be always clear from the context). We write $\si_j=(0,....,0,1,0,...)$ with the $j$\tss{th} element 1.
By convention, operations between sequences $\wek x=(x_j)$, $\wek y=(y_j)$ are almost always (see the exception $\sx^{\wek n}$ in the next subsection)  understood element-wise, or term-by-term, e.g., $
  \wek x\,\wek y=(x_j\,y_j),\quad \wek x+\wek y=(x_j+y_j)$ or $\wek x\vee\wek y=(x_j\vee y_j)$ (maximum), etc., whenever the component operations are well defined. By the same token, for the partial order $\wek x\le \wek y$ we require that $x_j\le y_j$ for every $j$.
  \vv

  A function is called {\em finitary}, if it takes only finitely many nonzero values. For example, a sequence $\sx$ is  finitary  if $x_j=0$ eventually.  We may also need the notion of 1-finitary sequences $\sx$, i.e., such that $x_j=1$ eventually. The set of finitary 0-1 sequences is denoted by $\BD$. For the sake of brevity we often suppress the zero tail, e.g., $101=101000...$. Conversely, finite sequences can be expanded or concatenated by additional 0s. The number of digits used in a finite sequence is our choice and does not necessarily coincide with the last 1. Therefore, the position's number  of the last digit in a finite sequence will be called its {\em depth} and the {\em depth} of a set (or a system) of sequences is understood as the maximum of depths of the set's members. To reiterate:  the depth is subject to choice, in contrast to the {\em length} by which we understand the position number of the last 1 in a 0-1 finitary sequence.
  As usual, the {\em support} of a sequence $\sx\in\R^{\N}$ is $\set{n\in\N:x_n\neq 0}$ and the support of a system of sequences is understood as the  supremum of all supports of the system's elements. \vv

  Addition modulo 2 is denoted by $p\oplus q,\,p,q\in\set{0,1}$.
  The {\em dyadic induction} is but the mathematical induction, tailored to fit the structure of $\BD$. That is, consider a property ${\cal P}(\dd)$. Then, if
\begin{enumerate}
\item ${\cal P}(\si_j)$ holds for every $j=1,2,...$,
\item for every $\sp,\sq\in\BD,\,\sp\sq=0$ the properties
${\cal P}(\sp)$ and ${\cal P}(\sq)$ imply the property ${\cal P}(\sp\oplus \sq)$,
\end{enumerate}
then ${\cal P}(\dd)$ is valid for every $\dd\in \BD$. In the second assumption it suffices to restrict to $\sp,\sq\in\set{\si_k,\, k=1,2,...}$.
\vv
We use the multiplication and summation brackets
\[
[\sx]\df\prod_n x_n,\qquad  \langle \sx\rangle \df\sum_n x_n.
\]
For the sake of clarity we will write the power $\wek x\supwek{n}=x_1^{n_1}\,x_2^{n_2}\cdots = [\wek x\supwek n]$  as an exception to the aforementioned convention
 (e.g., $\sx^{\son}=[\wek x]$). Then multivariate polynomials will appear as $\sum_{\wek n} \mbox{$a$}_{\wek n} \,\wek x\supwek n$. Occasionally (and very rarely) we may need the term-by-term powers, and then we use the following convention
 \be\label{brack}
 \{\wek x\}^{\wek n}=(x_1^{n_1}, x_2^{n_2},\dots).
 \ee
 Note that for a constant sequence $\wek x=(x,x,x,...)$ we have $\wek x^{\wek n}=x^{\langle \wek n\rangle}$.
By {\em polynomial chaos}, or ``chaos'' in short, we understand the sum of multilinear forms
\[
\sum_{\wek d\in \BD} a_{\wek d} \,\wek x\supwek {\sd} =\sum_{n=0}^\infty C_n, \quad
\mbox{with $n$-homogeneous chaos}\quad C_n=\sum_{\la\wek d\ra=n} a_{\wek d} \,\wek x\supwek {\sd}.
\]
The convention $0^0=1$ entails a number of relations (tautologies). For $\sp,\sd\in \BD$
\begin{equation}\label{basalg}
\sp^\sce\sp^\sd=\sp^{\sce+ \sd}= \sp^{\sce\vee \sd},\qquad (\son-\sp)^{\sd}=[\son-\sp \sd]=(\son-\sd)^{\sp}=\sz^{\sp\sd}.
\end{equation}
In particular,
\[
[1-\sd] =\sz^{\sd} \mbox{ (i.e., $\sp=\son$ above)},\qquad
\sp^{\sd}=[\son-\sd(\son-\sp)]=(\son-\sd)^{1-\sp}.
\]
Also, for a finitary $\sx$,
\[
\Ri{\sd\sx}=\sum_{\sp\in\BD}\sd^{\sp}\sx^{\sp}=(\son+\sx)^{\sd},\quad
\sum_{\sp\in\BD}\sd^{\son-\sp}\sx^{\sp}=[\sd+\sx],\quad [\sd+\sx]=\sx^{\son-\sd}(\son+{\sx})^{\sd}.
\]
\begin{remark}\rm
We may prove the middle statement by replacing temporarily $\sd:=\sd+\se$, with a sequence $\se>0$ of positive numbers, and later let $\se\to 0$. In other words, the ``division by 0'' in the formula `` $\sd^{\son-\sp}\sx^{\sp}= \sd (\sx/\sd)^{\sp}$\,'' is justified since the final outcome does not require the non-zero restriction.
\end{remark}
\subsection{Rademacher chaos, Walsh series, Riesz products}
The square wave, i.e., the periodic odd extension $r(t)$ of $1\hspace{-3pt} {\rm I}_{[0,1)}$, being dyadically condensed and then truncated, entails {\em Rademacher functions} on $[0,1]$:
\[
r_n(t)=r(2^nt)\Big|_{[0,1]}.
\]
(By convention, $r_0(t)=1$.)
Their ``teeth, broken off'' (usually suitably normalized) are called {\em Haar functions}, while their powers (or generalized products)  are called {\em Walsh functions},
\[
w_n=\sr^\sp, \quad\mbox{via the dyadic representation} \quad
n=\sum_{j=1}^\infty p_j \,2^{j-1},
\]
with the natural order on $\N$ coinciding with the lexicographic order on $\BD$. Walsh functions form an orthonormal basis of $L^2$, and also a basis of $L^p$, $1< p<\infty$.  The ``block convergence'' of Walsh series (or Rademacher chaos) follows immediately from the Doob's martingale convergence theorem, That is,
for every $f\in L^p[0,1]$, $1\le p<\infty$, a.s. and in $L^p$,
\be\label{Doob}
f=\sum_n a_n\,w_n=\sum_{\dd\in \BD} a_\dd\, \wek r^\dd\df \lim_n \ce{f}{\cD_n},
\ee
where $\cD_n$ is the $n$\tss{th} dyadic field, i.e., the set algebra spanned by $r_1,\dots,r_n$. In fact, the indicator of every dyadic interval $[j2^{-n},(j+1)2^{-n})$ is equal to the scaled Riesz product $R(\wek s)=2^{-n}[1+\wek s\wek r]$ for some unique sequence $\wek s$ of signs of length $n$. If $f$ has a finite Walsh series, i.e., $f$ is measurable with respect to some $\cD_n$, then we may choose such $n$ and call it the {\em depth} of $f$ (which is consistent with the depth of a 0-1 sequence).\vv

For $p>1$, the ``block convergence'' is equivalent to the summability of a Rademacher chaos, in particular, to the ``linear convergence'' of the Walsh series. In contrast (cf.\ \cite[Chap.\ 8]{Szu}), an $f\in L^1$ may have a divergent Walsh series. For $p<1$ even the zero function can be represented by a nontrivial Walsh series.\vv

Walsh functions span the algebra $\BL$ of finitary functions, where the average (mathematical expectation) yields the inner product
\be\label{innprod}
\E \,r^{\sp}=[1-\sp]\quad\imp \quad \la \sr^\sp,\sr^\sq\ra= [1-\sp\oplus \sq]
\ee
(recall that ``$\oplus''$ denotes addition modulo 2).
 {\em Riesz products} $[1+\wek x\wek r]$, where $\sx$ are finitary,   form a total (i.e., linearly dense) subset of $L^p$, $1\le p<\infty$. \vv

 \begin{remark}\rm
The term ``{\em Riesz product}'' is adopted from \cite{SWS} although in the language of \cite{Par} or \cite{Mey} the object would be rather called a ``{\em discrete exponential vector}''  or ``{\em discrete coherent vector}''.
\end{remark}

 Hence, a linear operator, a.k.a. ``chaos operator'', can be examined through its actions either on Walsh functions or on Riesz products although the algebra of Riesz products is slightly more complicated. Put $x=\tanh u,\,y=\tanh v$. The operation on real numbers carries over to Riesz products. Defining
\[
x\diamond y\df\tanh(u+v) =
 \left\{
\begin{array}{ll}
\frac{x+y}{1+xy}, & xy\neq -1\\
0 & xy=-1
\end{array}\right.,
\]
we have
\[
[\son+\sx\sr]\,[\son+\sy\sr]=
[\son+\sx\sy]\,[\son+\sx\diamond \sy\,\,\sr].
\]
\begin{lemma}\label{phi}
Consider a real function $\phi$ and a finite $F\subset\BL$. Then there exists a polynomial $\tilde{\phi}$ such that $\tilde{\phi}(f)=\phi(f)$ for each $f\in F$.
In particular, if $f\in\BL$ then $\phi(f)\in \BL$ and $T\phi(f)=T\tilde{\phi}(f)$ for every mapping $T:\BL\to\BL$.
\end{lemma}
\Proof
 W.l.o.g. we may and do assume that $\phi(0)=0$.  Choose $n$ to be greater than the maximum of the supports of $f$, $f\in F$. Recall Riesz product representations $R(\wek s)$ of dyadic intervals (beneath  \refrm{Doob}). The polynomial interpolation on  $\set{\phi(f\,R(\wek s)):f\in F, \wek s\in\set{\pm 1}^n}$ entails a polynomial $\tilde{\phi}$ such that
 \[
\tilde{\phi}(f)\cdot R(\wek s)= \tilde{\phi}(f\,R(\wek s))=\phi(f\,R(\wek s))=\phi(f)\cdot R(\wek s),
 \]
where we utilized the strengthened assumption on $\phi$. The remaining statements are immediate consequences of this observation.
\QED
\section{Simple chaos operators}
Although we trace the footsteps of P-A. Meyer, we depart from his original notation:
\[
x_A=\prod_{i\in A} x_i,
\]
which for ($\pm 1$)-sequences leads to set-theoretic operations, e.g., $x_Ax_B=x_{A\oplus B}$ with the symmetric difference of sets $A\oplus B$. Alternatively, the tensor notation could be utilized, e.g., $x^{\otimes A}$ for the above product. Although it is convenient for commutative $x$'s, it quickly becomes confusing when the elements of the underlying sequence, e.g., operators,  do not commute. Besides, the power notation (used by many authors, e.g., in \cite{Rud} as so called ``generalized powers'') allows a quick transparent formalism.
\subsection{Discrete chaos operators}
The simple yet crucial property of the Rademacher functions is that $r_j^2=1$ for every $j$.  Following P-A. Meyer \cite[Chap.II]{Mey} and  using the Parthasarathy's \cite{Par} terminology, we introduce the following chaos operators
\[
\begin{array}{rl}
\mbox{{\em\em number}} \mbox{ (or \em\em conservation}):    & {{N}}_j\,\sr^{\sp} =p_j \sr^{\sp},\\
&\rule{5pt}{0pt}\mbox{e.g., }  {{N}}_1\,r_1r_2=r_1r_2, \,{{N}}_3\,r_1r_2=0                       \\
\mbox{{\em\em annihilation}}:    & D_j\,\sr^{\sp}     =p_j\,r_j \sr^{\sp}=p_j\,\sr^{\sp+\si_j},\\
&\rule{5pt}{0pt}\mbox{e.g., }  D_1\,r_1r_2=r_2,  \,D_3\,r_1r_2=0       \\
\mbox{{\em\em creation}}    :    & ^1\!D_j\,\sr^{\sp} =(1-p_j)\,r_j\,\sr^{\sp}=(1-p_j)\,\sr^{\sp+\si_j},
\\
&\rule{5pt}{0pt}\mbox{e.g., }  {^1\!D}_1\,r_1r_2=0,\,     {^1\!D}_3\,r_1r_2=r_1r_2r_3.
\end{array}
\]
Additional operators will appear frequently:
\[
\begin{array}{rll}
 \mbox{``replacement''}
 & R_j=D_j+\,^1\!D_j:      & R_j\,\sr^{\sp}     = r_j \,\sr^{\sp}=\sr^{\sp+\si_j},\\
 &&\rule{5pt}{0pt}\mbox{e.g., }  R_1\,r_1r_2=r_2,       \, R_3\,r_1r_2=r_1r_2r_3           \\
  \mbox{``anti-number''} & ^1\!{N}_j=1-N_j:   & ^1\!N_j\,\sr^{\sp}  = (1-p_j) \,\sr^{\sp},\\
 &&\rule{5pt}{0pt}\mbox{e.g., }  ^1\!N_1\,r_1r_2=0, \, {^1\!N}_3\,r_1r_2=r_1\,r_2    \\
  \mbox{``symmetry''} & S_j=1-2N_j:  & S_j\,\sr^{\sp}        = (1-2p_j)  \,\sr^{\sp},\\
 &&\rule{5pt}{0pt}\mbox{e.g., }  S_1\,r_1r_2=-r_1\,r_2, \, S_3\,r_1\,r_2=r_1\,r_2.    \\
   \mbox{``asymmetry''} & A_j={^1\!D_j-D_j}:  & A_j\,\sr^{\sp}        = (1-2p_j)  \,\sr^{\sp+\si_j},\\
& &\rule{5pt}{0pt}\mbox{e.g., }  A_1\,r_1\,r_2=-r_2, \, A_3\,r_1\,r_2=r_1\,r_2\,r_3.    \\
 \end{array}
\]
The left superscript `1' toggles between two states of an operator: the original and the orthogonal complement for the number operator, or the original and adjoint for the annihilation operator (see the proposition below). The left superscript `0' leaves the operator unchanged.
\vv

Let us record redundantly convenient multiplication table (the subscript, or variable, is fixed and suppressed):
\begin{center}
\begin{tabular}{c|ccccc}
    {$_{\sf left}\!\backslash\!^{\sf right}$}\!    & ${N}$ & $R$   & $D$ & $^1\!D$   & $^1\!N$ \\ \hline
${N}$   & ${N}$ & $^1\!D$   & $0$   & $^1\!D$   & $0$     \\
$R$   & $D$ & $1$     & ${N}$   & $^1\!N$ & $^1\!D$   \\
$D$   & $D$ & $^1\!N$ & $0$   & $^1\!N$ & $0$     \\
$^1\!D$   & $0$   & ${N}$   & ${N}$ & $0$     & $^1\!D$   \\ 
$^1\!N$ & $0$   & $D$   & $D$ & $0$     & $^1\!N$ \\
\end{tabular}
\end{center}
\begin{proposition}
All seven operators are contractions with respect to the operator norm in every $L^p$, $1\le p\le\infty$. Furthermore,
\begin{enumerate}
\item
$R_j$, $N_j$, $^1\!N_j$, $S_j$ are self adjoint, $A_j$ is skew adjoint,  while $D_j$ and $^1\!D_j$ are  mutually adjoint;
\item
 Every two operators with distinct indices commute and, with exception of the pairs $(N,{^1\!N})$ and $(S,N)$  every two distinct operators with equal indices do not commute;
\item
Compositions $\sR^\sd,\,^\sce \sN^\dd=(^\sce \sN)^\dd,\,\wek S^\sd,\,^\sce\sD^\sd=(^\sce \sD)^\dd$ inherit the aforementioned properties. Note that $^\sce \sN\ort={^{1-\sce}\sN}$ and $^\sce \sD^*= {^{1-\sce}\sD}$. Also, ${^{\dd-\sce}\sN^\dd}={^{1-\sce}\sN^\dd}$ and ${^{\dd-\sce}\sD^\dd}={^{1-\sce}\sD^\dd}$.
\end{enumerate}
\end{proposition}
\Proof
For a fixed index $j$ consider the orthogonal decomposition
\[
\BL=\BL_j +\BL\ort_j ,
\]
where  $\BL_j={\rm lin}\set{\sr^\sq:q_j=0}$ ($r_j$ absent) and  $\BL\ort_j={\rm lin}\set{\sr^\sq:q_j=1}=r_j\BL_j$ ($r_j$ present). The number operator $N_j$ is the projection $\ce{\cdot}{\BL\ort_j}$ onto $\BL\ort_j$ (orthogonal in $L^2$) while $^1\!N_j= \ce{\cdot}{\BL_j}$ (the orthogonal complement $N\ort_j$ in $L^2$).
\vv
Clearly, $R_j$ and $S_j$ are isometries. The contractivity of the remaining operators follows from the compositions:
\be\label{RSN}
D_j=R_jN_j,\quad{^1\!D_j}=\, R_j\,{^1\!N_j}=N_jR_j,\quad A_j=R_jS_j.
\ee
The operators and their compositions can be defined either on the Walsh basis or on Riesz products, and then extended by linearity to finitary functions. So the adjointness can be verified on either basis.
For example, in proving that ${^1\!D_j}=D^*$ it suffices to consider only $\la D\,\sr^\sp,\sr^\sq\ra=\la \sr^p , {^1\!D}\,\sr^\sq\ra$ and apply the tautology  entailed by \refrm{innprod}.
The remaining properties follow a fortiori.
\QED
\vv

The multiplication table and the suitable mappings (the subscript is fixed and suppressed) show that the spanned algebra is isomorphic to $B(\BC^2)$:
\[
 b^- =\ddmat 0 1 0 0             \leftrightarrow D,  \quad
 b^+ = \ddmat 0 0 1 0            \leftrightarrow {^1\!D}, \quad
 n=b^o=\ddmat 0 0 0 1             \leftrightarrow {N}.
 \]
 Pauli spin matrices $\sigma_1,\,\sigma_2,\,\sigma_3$ constitute a basis of Hermitian 2$\times 2$ complex matrices:
\[
\ddmat {t+z}{ x-\imath y}{x+\imath y}{t-z}=t+x\,\sigma_1+y\,\sigma_2+z\,\sigma_3,
\]
so they express as follows:
\be\label{Pauli}
  \sigma_1=\ddmat 0 1 1 0             \leftrightarrow R, \qquad
  \sigma_2=\ddmat 0 {-\imath} \imath 0\leftrightarrow \imath\,A,\qquad
   \sigma_3=\ddmat 1 0 0 {-1}          \leftrightarrow S.
\ee
The pairs $\set{N,R},\,\set{D,R},\set{D,{^1\!D}}$, $\set{S,R}$ generate the 4-dimensional algebra $\cong B(\BC^2)$ while the pairs $\set{N,{^1\!D}},\,\set{N,D}$ generate a 2-dimensional subalgebra $\cong \mbox{\bf{$\mathfrak a\mathfrak f\mathfrak f$}}(1)$ (cf. \cite[1.1.Ex.10]{Sam}).\vv

\begin{remark} \rm
 For a time being we set aside the issue of convergence. We desire a norm stemming from an inner product such as the Hilbert-Schmidt norm, as opposed to the operator norm. However, none of our operators is Hilbert-Schmidt on $L^2$ unless we truncate the domain to $L^2(\ca D_n)$ (cf. \refrm{Doob}) and examine the entire injective system. We'll return to this issue in Section \ref{RHS}.
 \end{remark}\vv
 Let us record the actions of powers. Recall the exception from the ``term-by-term'' operation rule: the powers $\sx ^{\sd}$ are products, not sequences (cf.\ the paragraph before \refrm{brack}).  It suffices to consider only 0-1 powers, because all operators are either nilpotent of order 2, or idempotent, or square roots of the identity.  Below, $\sp$ and $\sd$ are 0-1 finitary sequences, $\sx$ is a sequence of scalars (or of elements of a commutative multiplication system). As usual, for a nonzero operator $T$, we write $T^0=I$.

\begin{proposition}
For $\wek p,\wek d\in \BD$, we have
\be\label{opWalsh}
\begin{array}{rlrl}
 {\wek{N}}^{\sd} \sr^{\sp}  &=\sp^{\sd} \sr^{\sp},
   &{\wek{N}}^{\sd} \sxr & =(\sx\sr)^{\sd} (1+\sx\sr)^{\son-\sd},\\
{^1\wek N} ^{\sd}\sr^{\sp} &=(\son-\sp)^{\sd}\sr^{\sp},
&{^1\wek N}^{\sd} \sxr&=(1+\sx\sr)^{\son-\sd}\\
 \sR^{\sd} \sr^{\sp} &=\sr^{\sp+\sd},
   &\sR^{\sd} \sxr & =\sr^{\sd}\sxr,\\
 \sD^{\sd} \sr^{\sp} &={\sp}^{\sd}\, \sr^{\sp+\sd},
   &\sD^{\sd} \sxr&=\sx^{\sd} (1+\sx\sr)^{\son-\sd},\\
^1\sD^{\sd} \sr^{\sp} &=(\son-\sp)^{\sd} \,\sr^{\sp +\sd},
   &^1\sD^{\sd} \sxr&=\sr^{\sd} (1+\sx\sr)^{\son-\sd}.\\
   \wek S^{\sd} \sr^{\sp} &=(\son-2\sp)^{\sd} \,\sr^{\sp},
   &\wek S^{\sd} \sxr&= [1+(-1)^\dd\sx\sr].\\
\end{array}
\ee
All powers at once can be captured by the operator Riesz products yielding more transparent formulas:
\be\label{opRiesz}
\begin{array}{c}
\vspandexsmall\Ri{\su\sN} \,\sxr = \Ri{(\su+1)\sx\sr},\quad
\Ri{\su\,^1\sN}\,\sxr  =[1+\su+\sx\sr],\\
\vspandexsmall\Ri{\su\sR}\,\sxr  =\Ri{\su\sr}\,\sxr,\\
\vspandexsmall\Ri{\su\sD}\,\sxr  =[1+\su\sx+\sx\sr],\quad
\Ri{\su\,{^1\sD}}\,\sxr  =\Ri{(\su+\sx)\sr},\\
\vspandexsmall\Ri{\su\wek S}\,\sxr  =\Frac{\big([1+\su]+[1-\su]\big)\,[1+\sx\sr]+
\big([1+\su]-[1-\su]\big)\,[1-\sx\sr]}{2}.
\end{array}
\ee
\end{proposition}
\Proof
 The formulas  follow directly from the definitions of the operators and from 0-1 algebra.
 \QED
 \vv

The Ornstein-Uhlenbeck semigroup of operators, bounded in $L^p$ for $1\le p<\infty$,  provides an example of chaos operators with each $n$-homogenous chaos being an eigenspace:
\be\label{hyp}
T_t: \sxr \mapsto \Ri{e^{-t}\,\sx\sr}\quad \mbox{ or }\quad \sr^{\sp}\mapsto e^{-t\la\sp\ra}\,\sr^{\sp}.
\ee
It entails  the unbounded generator $A=-\la \sN\ra$, akin to the Ornstein-Uhlenbeck generator \refrm{Mall}:
\[
 A=\lim_{t\searrow 0} \frac{1}{t}(T_t-I),\quad
 \sxr \mapsto -\sxr\,\sum_j \frac{x_j\,r_j}{1+x_j\,r_j}\quad\mbox{or}\quad
 \sr^{\sp}\mapsto -\la\sp\ra\,\sr^{\sp}.
 \]
If $\sr$ were a continuous variable, then $\sD^\dd$ would be exactly the partial derivative of order $\dd$ with respect to $\sr$. Further, our $\sD$ plays the role of the gradient, tied to ``the directional'' derivative in the ``direction'' $\sx$. Bringing temporarily the usual dot product notation, we observe that
\[
      \sx\cdot \sD=\la \sx\sD\ra=\sum_j x_j D_j=\frac{d}{ds}\Ri{s\,\sx\sD}\Big|_{s=0}.
\]
The ``direction'' could be a sequence of operators, yielding yet another chaos operator through the ``directional derivative''. For example, the sequence $x_j=c_jR_j$, where $c_j$ are scalars, defines the self-adjoint operator
      \[
      D_{\sce}\df \sce\sR\cdot \sD=\la\sce\sN\ra=\sum_j c_j N_j =\sce\cdot \sN,
      \]
and thus $ D_{\sce} \sr^{\sp}=\la\sce\sp\ra\, \sr^{\sp}$. The action on Riesz products yields an unwieldy formula but in the differential form, since $\sN=\sR\sD$, it is transparent due to \refrm{opRiesz}:
\[
D_{\sce}[1+\sx\sr]=\frac{d}{ds} [1+s\,\sce\sN]\,[1+\sx\sr]\,\Big|_{s=0}= \frac{d}{ds} [1+(s\,\sce+1)\sx\sr]\,\Big|_{s=0}.
\]
\begin{example}{~}

 \begin{enumerate}
        \item Let ``$\oplus$'' denote addition modulo 1 on $[0,1]$.  The {\em dyadic derivative} {\rm \cite[p.40]{SWS}} $D=D_{\sce}$ is obtained when $c_j=2^{j}$. Its direct action follows the formula:
            \[
            D f(x)=\frac{1}{2}\sum_{j=1}^\infty\frac{f(x)-f(x\oplus 2^{-j})}{2^{-j}}.
            \]
        Under the lexicographic (Paley's) order, in the language of Walsh-Fourier series
        \[
        f=\sum_n a_n w_n\quad \mapsto \quad D^k f=\sum_n n^k a_n w_n,\quad \mbox{i.e., } \widehat{D_k f}_n =n^k \widehat{f}_n,
        \]
        as one would expect from the Fourier transform of the $k^{\rm th}$ derivative.
        \item Let $\sce=\son$, i.e. $D_{\son}=\sum_j N_j$. Then
        \[
        D_{\son}^k\, r^{\sp}=\la\sp\ra^k \, r^{\sp}\quad\mapsto\quad e^{-tD_1} =T_{t},
        \]
        which is the hypercontraction operator \refrm{hyp}.
      \end{enumerate}
 \end{example}
\subsection{Leibnitz formula and Chain Rule}
We will derive complete analogs of the Leibnitz formula for the ``differential operators''  $\sN$ and $\sD$ (preliminary and limited versions were obtained in \cite[Chap.\ 12.2]{Szu}).  Recall the symmetry $S_j={^1\!N_j}-N_j=1-2N_j$ and the asymmetry $A_j={^1\!D}_j-D_j$.
\begin{lemma} \label{lem:Nj}
Let $ f,g\in \BL$. Fix $j$ and suppress the subscript (e.g., $N=N_j$).
\begin{enumerate}
\item \label{Nj}
$   N (fg)=\big(N f\big)\,\big({^1\!N} g\big)+\big({^1\!N}f\big)\,\big(N g\big)$,

$   ^1\!N (fg)=\big({^1\!N} f\big)\,\big({^1\!N} g\big)+\big({N}f\big)\,\big(N g\big)$;
\item
Hence, $S (fg)=(Sf)(Sg)$. Therefore, $S(f^n)=(Sf)^n$. Thus
\[
N f^n=\Frac{f^n-(Sf)^n}{2},\quad {^1\!N} f^n=\Frac{f^n+(Sf)^n}{2},
\]
or capturing both states of the number operator at once,
\[
{^e\!N}\, f^n=\Frac{f^n+(2e-1)(Sf)^n}{2}
\]
\item
For every real function $\phi$, $S(\phi(f))=\phi(Sf)$ and the ``chain rule'' holds:
\[
N \phi(f)=\frac{\phi(f)-\phi(Sf)}{2}\quad\mbox{and} \quad
\,{^e\!N} \phi(f)=\frac{\phi(f)+(2e-1)\phi(Sf)}{2}.
\]
\item
$   N (fg)=R\left[\big(D f\big)\,\big({^1\!D} g\big)+\big({^1\!D}f\big)\,\big(D g\big)\right]$,

$   ^1\!D (fg)=\left[\big({^1\!D} f\big)\,\big({^1\!D} g\big)+\big({D}f\big)\,\big(D g\big)\right]R$;
\end{enumerate}
\end{lemma}
 \Proof
The first formula follows from the relation $p\oplus q=p(1-q)+(1-p)q$ on $\set{0,1}$, which by definition of the operators $N_j$ and ${^1\!N}$, acting on $\set{\sr^{\sp}}$, entails the identity
\[
N\big(\sr^{\sp} \sr^{\sq}\big)
   =\big(N \sr^{\sp}\big)\,\big({^1\!N}\sr^{\sq}\big)+ \big({^1\!N}\sr^{\sp}\big) \big(N\sr^{\sq}\big).
\]
Hence, with finite Walsh series, representing functions $f$ and $g$,  by linearity, we arrive at the full identity:
   \[
   \begin{array}{c}
   N(f\cdot g)=N\left(\Sum_{\sp} a_{\sp} \sr^{\sp}\cdot \Sum_{\sq} b_{\sq}\sr^{\sq}\right)=\Sum_{\sp}\Sum_{\sq} a_{\sp}b_{\sq} N\big(\sr^{\sp} \sr^{\sq}\big)\\
  =\Sum_{\sp}\Sum_{\sq} a_{\sp}b_{\sq} \big(N \sr^{\sp}\big)\,\big({^1\!N}\sr^{\sq}\big)
  +\Sum_{\sp}\Sum_{\sq} a_{\sp}b_{\sq} \big({^1\!N}\sr^{\sp}\big) \big(N\sr^{\sq}\big)\\
  =\big(N f\big)\,\big({^1\!N} g\big)+\big({^1\!N}f\big)\,\big(N g\big).
  \end{array}
   \]
The dual formula for $^1N$ follows by the same token.\vv

In Statement (2), using the relations $S=1-2N=2{^1\!Nf}-1$, we obtain the formula
\be\label{Nf2}
N(f^2)=2\big(Nf\big)\big({^1\!Nf}\big) =\frac{1}{2} (f-Sf)(f+Sf)=\frac{f^2-(Sf)^2}{2}.
\ee
This, in turn implies that $S(f^2)=(Sf)^2$, which by polarization yields $S(fg)=(Sf)(Sg)$ and, consequently, $S(f^n)=(Sf)^n$. In other words,  $S$ and the power operation commute, i.e., $S$ is an isomorphism. Now, the formula for $N(f^n)$ follows by simple induction. The twin formula may be repeated or derived directly from the relation between two states of the number operator.\vv

For a general $\phi$ in statement (3) we recall Lemma \ref{phi}.
\vv
The formula for $D(fg)$ in Statement 4 holds up to the factor $R$, since $D=RN$ and ${^1\!D}=R{^1\!N}$.
\QED

\begin{theorem}\label{th:sN}
Let $ f,g\in \BL$, $\dd,\,\se\in\BD$.
Then
\begin{enumerate}
\item
${\sN}^\dd(fg)=
\Sum_{\sce\le \sd} \Big(  {^{\sce}{\sN}^{\dd}} f  \Big)\,\Big({^{\dd-\sce}{\sN}^\dd} g\Big)$,

$\,{^\se\sN}^\dd(fg)=
\Sum_{\sce\le \sd} \Big(  {^{\sce}{\sN}^{\dd}} f  \Big)\,\Big({^{\dd-\se\oplus \sce}{\sN}^\dd} g\Big)$;
\item
${\sD}^\dd(fg)=
\sR^\dd\Sum_{\sce\le \sd} \Big(  {^{\sce}{\sD}^{\dd}} f  \Big)\,\Big({^{\dd-\sce}{\sD}^\dd} g\Big)$,

$\,{^\se\sD}^\dd(fg)=\sR^\dd
\Sum_{\sce\le \sd} \Big(  {^{\sce}{\sD}^{\dd}} f  \Big)\,\Big({^{\dd-\se\oplus \sce}{\sD}^\dd} g\Big)$.
\end{enumerate}
\end{theorem}
\Proof
The tautology
\[
\sp\dd=0,\,\sq\le\sp,\,\sce\le\dd\quad\imp\quad ^\sq\sN^\sp\,^\sce\sN^\dd=\,^{\sq+\sce}\sN^{\sp+\dd}.
\]
 entails the dyadic induction. That is, assuming that the formula is valid for $\dd$ and  $\sp$ such that $\sp\sd=0$, we apply $\sN^\sp$ to both sides of the equality.  By the assumption  the right hand side equals
\[
\Sum_{\sce\le \sd,\sq\le\sp}\left( {^\sq\sN^\sp}\, {^{\sce}{\sN}^{\dd}} f  \right)\left({^{\sp-\sq}{\sN}^\sp}\, {^{\dd-\sce}{\sN}^\dd} g\right)
=
\Sum_{\sq+\sce\le \sp+\sd} \left( {^{\sq+\sce}\sN^{\sp+\sd}}  f  \right)\left({^{\sp+\dd-(\sq+\sce)}{\sN}^{\sp+\sd}}  g\right),
\]
and since for every $\wek e\le \sp+\dd$ there are unique $\sq\le\sp$ and $\sce\le\dd$ such that $\wek e=\sq+\sce$, the formula holds for $\sN^{\sd+\sp}$, as required. \vv

Consider now a single variable $j$ and observe in Lemma \ref{lem:Nj}.\ref{Nj} that the left superscripts $c=u_j$ and  $c=1-u_j$ appear in the first factors ${^cN f}$. Then the left superscript in the second factor in each summand is $1-c\oplus e$. The twin formula repeats the pattern in the $\dd$-based composition.\vv

The remaining formulas follow from the compositions $D=RN$ and ${^1\!D}=R\,{^1\!N}$.\QED
\begin{corollary}\label{cor:eNphif}
Let $\wek e, \dd\in\BD$, $f\in\BL$, and $\wek u$ be a finitary sequence of real numbers. Then
\begin{enumerate}
\item
for a real function $\phi$,

$
{^\wek e}\sN^\dd \phi(f)=2^{-\dd}\sum_{\sce\le\dd}(-1)^{\wek c(1-\wek e)}\phi(\wek S^\sce f)$, and for a Riesz product,

$\left[1+2\wek u \,{^\se \sN}\right]\phi(f)=\left[1+\wek u\big(1+(2\se-1)\phi(\wek S f)\big)\right]$;
\item
for a  sequence $\wek f\in \BL^{\N}$ and a finitary sequence $\sp$ of integers ($p_k$ may be real if $f_k>0$),

${^\se\sN}^\dd\wek f\,^\sp=\left(\Frac{1-(1-2\se)\{\wek S \wek f\}^\sp}{2}  \right)^\dd$, and for a Riesz product,

$\left[1+2\wek u \,{^\se \sN}\right]\wek f\,^\sp=\left[1+\wek u\left(1+(2\se-1)\{\wek S \wek f\}^\sp\right)\right]$.
\end{enumerate}
\end{corollary}
\Proof
 First, we use ${^\wek e}\sN^\dd=\left(\Frac{1+(2\se-1)\wek S}{2}\right)^\dd$ and the identity  $(2\wek e-1)^\sce=(-1)^{\wek c(1-\wek e)}$. Then we use the commutativity of $\wek S$ with every product together with Lemma \ref{phi}. \QED
 \vv
The analogs of formulas in Corollary \ref{cor:eNphif} for the creation and annihilation operators, in terms of the symmetry $\wek S$ can be obtained simply by multiplication by $\wek R^\dd$, as done in Theorem \ref{th:sN}. However, if we wish to refer to the asymmetry  $\wek A={^1\sD}-\sD=\wek R\wek S$, we observe that it does not commute with powers (or with $\phi$). Indeed, using the decomposition $\phi=\phi_e+\phi_o$ into the sum of even and odd parts (e.g., $e^t=\cosh t+\sinh t$), we see that
\[
A\phi(f)=\phi_e(Af)+R\phi_o(Af).
\]
Following this approach, we further check that
\[
\phi(Rg)=\phi_e(g)+R\phi_o(g)\quad\imp\quad \phi(\wek A^\sce f)=\phi_e(\wek A^\sce f)+\sR^\sce\phi_o(\wek A^\sce f).
\]
This change will affect the power formulas, e.g.,
\[
 {^\se\sD}^\dd\phi(f)=\left(
 \Frac{1+(2\se-1)\Big(\phi_e(\wek A f)+\wek R \phi_o(\wek A f)\Big)}{2}
 \right)^\dd
 \]
 or for a Riesz product,
 \[
 \left[1+2\wek u \,{^\se \sD}\right]\phi(f)=\left[1+\wek u\Big(1+(2\se-1)\Big(\phi_e(\wek A f)+\wek R \phi_o(\wek A f)\Big)\Big)\right].
 \]
 It is not clear whether this angle of view brings much of value, which suggests to stay with simple and transparent expressions based on the symmetry $\wek S$, or to confine to the class of either even or odd functions. Nevertheless, the switch from $S$ to $A$ can be performed at will.
\section{Signed groups}
The next natural step is to investigate the algebra spanned by compound operators such as
${N}_1R_2{^1\! D}I_3N_4D_5 A_6S_7\cdots$.
\vv
Each of operators $\sN, \, \sR,\,\wek S$, together with $\pm I$,  generates a commutative group, then the algebra, that is isomorphic to the group of Walsh functions $\set{\pm\sr^\sp}$ generated by the Rademacher functions. So, it is natural to call each of them a {\em Rademacher system}. In contrast, $\ca C\df\set{\pm\sR^{\wek s}\,{\wek S}^\sp: \wek s,\wek p\in\BD}$ becomes a group of mixed commutativity with the action
\be\label{action}
\sR^{\wek t}\wek S^{\wek p}=(-1)^{\wek t\wek p}\wek S^{\wek p}\sR^{\wek t}
\quad\imp\quad
\Big(\sR^{\wek s}\wek S^\sp\Big)\,\Big(\sR^{\wek t}\wek S^\sq\Big)=(-1)^{\wek p\wek t}\,\sR^{\wek s+\wek t}\wek S^{\sp+\sq}
\ee
(recall that for a constant sequence $a^\dd=a^{\la\dd\ra}$). In particular, $\left(\sR^{\wek t}\wek S^\sp\right)^2=(-1)^{\wek t\sp}\in\set{-1,1}$. If we confine to $\BD_n=\set{-1,1}^n$, then the corresponding group $\ca C_n$ of order $2^{n+1}$  is a basis of the spanned algebra $\cC_n$. Therefore, $\ca C$ is a Hamel basis of the algebra $\cC$ of finite rank linear operators on the vector space  $\BL$ of functions on $[0,1]$ with dyadic supports. The motivation behind the choice of the acronym is clear: ``C'' for ``chaos''.
\vv
Extracting the few crucial properties, we call an associative multiplicative structure $\ca S$ {\em signed} if
\begin{enumerate}
\item all members either commute or {\em anticommute} (i.e., $ss'=-s's$);
\item $\ca S$ contains the group $\set{-1,1}$ whose members are called {\em  scalars}, and $1e=e1=e\neq (-1)e=e(-1)$ (denoted simply by $-e$);
\item For every member $e$,  $e^2\in\set{-1,1}$.
\end{enumerate}
Thus, elements of $\ca S$ can be called either {\em  positive} or {\em  negative}, depending on the sign of their squares.
The adjectives
``signed'', ``positive'' or ``negative'',``commutative'' or ``anticommutative'',  extend to the sets whose all members
have the corresponding property. By convention,  a singleton is both commutative and anticommutative.   If all signs of members of a signed set $E$ are
equal, then we call $E$ {\em  pure}, otherwise we call $E$ {\em mixed}. For example (cf. \refrm{Pauli} or Example \ref{pq} below), the set $\set{\sigma_1,\sigma_2,\sigma_3}$ of Pauli matrices  or $\set{R,\imath A,S}$, being all positive, and the set $\set{i,j,k}$ or $\set{\imath S,-A,\imath R}$ of quaternions, being all negative, are pure. However, $\set{R,S,A}$ is mixed because $R, S$ are positive and $A=RS$ is negative. \vv

The signum function $\sigma(e)=e^2\in\set{-1,1}$ is well defined on $\ca S$ as well is the symmetric commutativity function on $\ca S\times \ca S$:
 \[
c(e,e')=e\,\circ\, e'=\left\{
\begin{array}{rl}
1,&\mbox{if~ $e\,e'=e'\,e$},\\
-1,&\mbox{if~ $e\,e'=-e'\,e$}.\\
\end{array}
\right.
\]
\begin{example}\label{commD}
In the group $\ca C$, $\left(\sR^{\wek s}\wek S^\sp\right)\,\circ\,\left(\sR^{\wek t}\wek S^\sq\right)=(-1)^{\wek p\wek t+\wek s\wek q}$.
\end{example}
Recall that the position $n$ of the last element of a 1-finitary sequence $\wek e$ that is different from 1 is called  its {\em length} (akin to yet distinct from ``depth'') and denoted by $\ell(\se)$.
\vv

By convention, $e^0=1$. A signed sequence $\wek e=(e_1,e_2,\dots)$ generates the signed group $ \ca G(\wek
e)=\set{\pm\,\wek e^{\wek p}:\sp\in\BD}$. We say that a signed set $E$ or sequence $\wek e$ is {\em  nonbasic} if for
some $k\ge 2$ there exist distinct members $e_1,\dots,e_k$ such that $e_1\dots e_k$ is a
scalar. Otherwise, we call the set or sequence {\em  basic}. In particular,
every singleton is basic and if a basic set contains at least two elements then it is scalar-free. Equivalently, $\wek e$ is basic if and only if for
every finitary $\wek p,\,\wek q\in \BD$, $\wek e ^{\wek p}=\wek e^{\wek q}$
implies $\wek p=\wek q$. Also, if $\ell(\wek e)\ge 3$, then $\wek e$ is nonbasic if and only if one element expresses as a product of two distinct elements.
We say that an element $x$ is {\em independent} of a set $E\leftrightarrow\wek e$, if no power $\wek e^{\wek p}=x$.\vv

Clearly, a basic generator of length $n$ entails the group of order $2^{n+1}$. Conversely, every nontrivial finite signed group in $\ca S$ must have order $2^{n+1}$ for some $n\ge 0$. The case $n=0$ corresponds to $\ca G=\set{\pm 1}$. From now on by ``generator'' we will understand ``basic generator''.
\subsection{Double logic}
Any doubleton $D=\set{a,b}$ supports the unary action $x\mapsto \mbox{not }x$. Let us choose $D=\set{-1,1}$.  There are eight one-to-one functions on the double doubleton $D^2$. Just two functions, with the notation borrowed from the corresponding Rademacher chaos operators,
\[
\begin{array}{ll}
(s_1,s_2) \stackrel{R}{\to} (s_2,s_1)     &\mbox{ (for ``Replacement'')},\\
(s_1,s_2) \stackrel{S}{\to} (s_1,-s_2)   &\mbox{ (for ``Sign Swap'')},
\end{array}
\]
entail all
\[
\pm\, R^sS^p,\quad s,p\in \set{0,1}.
\]
Equivalently, we may consider the symmetries of the square $(\pm 1,\pm 1)$: $R$ being the revolution in $\R^3$ (or symmetry in $\R^2$) about the main diagonal, and $S$ being the symmetry (or revolution in $\R^3$) about the horizontal axis.\vv

 Let us use $2\times 1$ vectors to indicate the operations: the lower cell codes the potential replacement and the upper cell codes the potential  sign swap:
\[
\dlog 0 1 = R,     \qquad \dlog 1 0=S.
\]
The multiplication tables follow, extending immediately to the Cartesian product functions on $\left(D^2\right)^{\N}$,
\begin{equation}\label{product}
\begin{array}{c}
\dlog s p\,\dlog t q =(-1)^{pt} \dlog {s\oplus t}{p\oplus q}\quad\imp\\
\dlog {\wek s} {\wek p}\,\dlog {\wek t} {\wek q} =(-1)^{\wek p \wek t} \dlog
{\wek s\oplus \wek t}{\wek p\oplus \wek q}, \quad
\dlog{\wek s}{\sp},\dlog{\wek t}{\sq}\in\BD^2
\end{array}
\end{equation}
(with addition modulo 2). Let us denote the group of such finitary double sequences by $\cal D$. In view of \refrm{action} it is isomorphic to the group $\ca C$  induced by the chaos operators $R$ and $S$. While the definition itself does not depend on the order of the sequence, in practice one must choose either of two orders. For reasons apparent in the forthcoming matrix setting, we enumerate sequences from right to left. In particular, a finite 0-1 sequence will by concatenated to the left by augmented 0s, as illustrated by the following examples. The added zeros correspond to the identity operations on the
third and forth coordinate (counting from right).
\[
\begin{array}{l}
\displaystyle\vspandex\dlog{1101}{1011}\dlog{11}{01}
=\dlog{1101}{1011}\dlog{{\textcolor{red}{00}}11}{\textcolor{red}{00}\,01}=\dlog{1110}{1010},\\
\displaystyle\dlog{01}{10}\dlog{1101}{0101}=\dlog{\textcolor{red}{00}01}{\textcolor{red}{00}11}\dlog{1101}{0101}=-\dlog{1100}{0110}.
\end{array}
\]
(Writing the first factor above the second factor helps to establish the sign).
 Clearly, if the supports of a double sequence $x=\dlog {\wek s} {\wek p}$ and of a set $E$ of double sequences are disjoint, then $x$ is independent of $E$ (but not conversely).\vv

 Formula \refrm{product} entails the signum function
$
\sigma \dlog {\wek s}{\wek  p} =(-1)^{ \wek s\,\wek p}$ and the commutativity function $\dlog {\wek
s}{\wek p}\mbox{\footnotesize $\circ$}\dlog {\wek t} {\wek q}=(-1)^{\wek p\,\wek
t+\wek q\,\wek s}$. Hence the group $\ca D$ contains infinitely many positive
and infinitely many negative elements.  \vv

Below we discuss the preservation of
the anticommutativity of a finite sequence  $\wek e=\left(\dlog{\wek s}{\sp}\right)$ while increasing the depths of its members and its length.
That is, $\dlog{\wek s}{\sp}\mapsto \dlog{s'\wek s}{p'\sp}$, with the same new column for all members. Then we add one more element, choosing the simplest $\dlog{s''\wek 0}{p''\wek 0}$ by selecting a new column on the left.
The new sequence of the length increased by 1 is denoted by $\wek e'$.
\begin{enumerate}
\item
Any choice of two distinct columns $\dlog 1 0,\,\dlog 0 1,\,\dlog 1 1$
preserves the anticommutativity.
\item Let $\wek e$ be negative. To keep $\wek e'$ negative the only choices are
$\dlog{1\wek 0}{1\wek 0},\,\dlog{1\wek s}{0\wek p}$ or $\dlog{1\wek 0}{1\wek 0},\,\dlog{0\wek s}{1\wek p}$.
\item
If $\wek e$ is negative, then in order to obtain a positive $\wek e'$ the only choices
are
$\dlog{1\wek 0}{0\wek 0},\,\dlog{1\wek s}{1\wek p}$ or $\dlog{0\wek 0}{1\wek 0},\,\dlog{1\wek s}{1\wek p}$
\item
If $\wek e$ is mixed or positive, it is not possible to switch to a
negative sequence $\wek e'$ by a single augmentation. A double augmentation is needed, as shown by the
example
$\dlog{10\wek 0}{10\wek 0},\,\dlog{01\wek s}{11\wek p}$.
\end{enumerate}
\vv

The number $N=N(\wek e)$  of negative elements in a basic sequence $(e_1,\dots,e_n)$, serving as a generator of a group, is not a group property, that is, it is not invariant under a group isomorphism. In contrast, the
total number of negative elements in the generated group and the commutativity matrix up to permutation of the corresponding rows and columns are group properties.

\begin{theorem}\label{embed}
Let $\wek e=(e_k)$ be a countable signed sequence in a multiplicative system
$M$. Let $\bm\sigma =\sigma(\wek e)$ and $c_{jk}=e_j\circ e_k$. Then the
group $\ca D$ contains a sequence $d_k=\dlog{\wek s_k}{\wek p_k} $
with the matching signs and the commutativity function.
\end{theorem}
\Proof Recall that $\dlog 01$ and $\dlog 10$ are positive while $\dlog 11$ is negative, and all three elements in this triple $T$ anticommute. For $d_1$ choose one from $T$ to match the sign $\sigma_1$. Then, choose $d_{21}\in T$ such that $d_{21}\circ d_1=c_{21}$, and choose $d_{22}\in T$ to define $d_2= \left(d_{22}\,d_{21}\right]$ with $\sigma(d_2)=\sigma_2$. Note that the added element $d_{22}$ does not affect commutativity.
The remainder of the construction follows by induction. That is, if $d_1,\dots, d_n$ have been constructed that match all the required signs, and $d_k$ has depth $k$ for $k=1,\dots, n$, then $d_{n+1}=\left(d_{n+1,n+1}\,\cdots d_{n+1, 1}\right]$, where the consecutive elements from right to left are selected from the triple $T$ to match the consecutive signs $c_{n+1,1},\cdots c_{n+1,1}$ and the last element $d_{n+1,n+1}$ makes the match with $\sigma_{n+1}$.
\QED

\begin{example}{~}\label{examples:embed}\rm
The embedding algorithm of Theorem \ref{embed} is flexible and can be modified in so many ways to accommodate desired patterns. As in its proof, we will be sampling the triple $T=\set{\dlog 01,\,\dlog 10,\,\dlog 11}$. Consider the following examples.
\begin{enumerate}
\item
A positive or negative anticommutative basic sequence:
\[
d_1=\dlog{1}{0},\,e_1=\dlog{1}{1},\quad d_{k}=\underbrace{\dlog{10\dots 0}{01\dots 1}}_{k},\quad
e_{k}=\underbrace{\dlog{10\dots 0}{11\dots 1}}_{k},\quad k\ge 2.
\]
\item
Infinitely many anticommutative sequences, mutually commutative:
\vv

Partition $\N$ into the union of infinitely many disjoint subsets, and on each subset perform the above construction, choosing desired signs.
\item
Two commutative sequences, mutually anticommutative:
\[
e_k=\dlog{\wek s_k 1}{\wek p_k 1},\quad e'_k=\dlog{\wek s'_k 1}{\wek p'_k 0},
\]
where $\wek s_k$ has the single 1 on the position $2k$, $\wek s'_k$ has the single 1 on the position $2k+1$. The bottom sequences also have a single position filled with 1 or 0 beneath the corresponding 1 in $\wek s$, and the  choice follows our wish of having both sequences positive (so we choose 1 for $\wek e$ and 0 for $\wek e'$) or negative (choosing $0$ for $\wek e$ and 1 for $\wek e'$).
\item
Finitely or infinitely many commutative sequences, mutually anticommutative:

\vv
First, we choose an anticommutative sequence $\wek c$, as in the first example above, but stretch its elements if necessary, keeping infinitely many (say, even) positions, marked as natural numbers $\N$, ready to accept modifications. Then we partition $\N$ into the union of disjoint infinite subsets arithmetic progressions. Over each subset we define a sequence $\wek d_k$ using a single choice from the triple $T$ on disjoint positions to ensure commutativity. The specific selection depends on the requested signs, and can be handled like above. Finally, we combine the disjoint $\wek d_k$ and $c_k$ from $\wek c$ to obtain the desired sequence $\wek e_k$.
\end{enumerate}
\end{example}

 \subsection{Matrix representation}
As an equivalent alternative to the square's framework from the previous subsection, we may place the vertices of our square at the coordinate axes so the operators become simply the matrices
 \[
 R=\ddmat 0 1 1 0, \quad S=\ddmat 1 0 0{-1}.
 \]
 The concatenation corresponds to the doubling or rather quadrupling of matrices. We denote the universal identity matrix by 1 and the universal zero matrix by
0. We can quadruple the identity matrix in four essential ways, of which the first one will be called {\em principal}:
 \def\hh{\rule{4pt}{0pt}}
 \def\hhh{\rule{5pt}{0pt}}
\[
\begin{array}{rlrl}
Q_0^0=Q_0^0(1)= &\begin{pmatrix}
\hh 1\hh&\hh0\hh\\
0&1\\
\end{pmatrix},\quad
Q_1^0=Q_1^0(1)= &
\begin{pmatrix}
\hh 0\hh &\hh 1\hh\\
1&0\\
\end{pmatrix}, \vspandex\\
\\
\vspandex
 Q_0^1=Q_0^1(1)=&
\begin{pmatrix}
\hh 1\hh &\hh 0\hh\\
0&-1\\
\end{pmatrix},\quad
Q_1^1=Q_1^1(1)=&
\begin{pmatrix}
\hh 0\hh &\hh 1\hh\\
-1&0\\
\end{pmatrix}.
\end{array}
\]
Then, for a quadratic matrix $e$ we denote
 \be\label{Qact}
Q_0^0(e)=\begin{pmatrix}
e&0\\
0&e\\
\end{pmatrix},\quad Q^s_p(e)=Q^s_p\,\cdot\, Q_0^0(e),
 \ee
with the actual matrix multiplication in the latter formula.
The iteration follows:
\be\label{action1}
Q^{\wek s}_{\wek p}=Q^{s_k,s_{k-1},\dots s_1}_{p_k,p_{k-1},\dots
p_1}=Q^{s_k}_{p_k}\Big(Q^{s_{k-1},\dots s_1}_{p_{k-1},\dots p_1}   \Big),
\ee
entailing $2^k\times 2^k$ matrices, quite scarce, of various depths $k$. Let us denote by $\ca M$ the group of such matrices spanning the algebra $\cM$.  When finitely (or even infinitely) many of these matrices interact, they can be brought to the same size by repeated doubling $Q^{...00}_{...00}$, corresponding to the previous catenation to the left.
Further, mere indices would suffice, for they satisfy the multiplication tables \refrm{product}. Dropping the symbol ``$D$''  induces the isomorphism with the group $\cal D$ (or the one generated by \refrm{action}).
\begin{remark}
{\rm From now on we will allow the convenient ambiguity by writing $\dlog{\wek s}{\wek p}$ (where sequences are ordered from right to left) in lieu of $Q^{\wek s}_{\wek p}$ (where sequences are ordered from left to right). The exact meaning will be clear from the context.}
\end{remark}
\begin{example}\label{pq}
\rm
Consider depth 2, i.e. $4\times 4$ matrices. Put $\imath=\dlog{1}{1}$ ($=\dlog{\textcolor{red}{0}1} {\textcolor{red}{0}1}$, using concatenation).   The Pauli matrices \refrm{Pauli}  in our notation appear as follows:
\[
\sigma_1=\dlog{00}{10},\, \sigma_2=-\dlog {11}{11},\, \sigma_3=\dlog {10}{00}.
\]
Then Pauli matrices, multiplied by $\imath$, yield the following basis of quaternions
 (cf. \cite{Zha}):
 \[
 \ddmat{t+\imath x}{y+\imath z}{-y+\imath z}{t-\imath x}=t+x\,{\bm i}+y\,{\bm j}+z\,{\bm k}.
 \]
 In our notation
\[
{\bm i}=\imath \sigma_3=\dlog {11}{01} , \,  {\bm j}=\imath \sigma_2=\dlog {10}{10}, \, {\bm k}=\imath\sigma_1=\dlog {01}{11} .
 \]
Among sixteen double sequences of depth 2,  $\imath$ commutes only with these six matrices, besides the identity and itself.
\end{example}

Denote by ${{c_{11}}}$ the number of occurrences of the element $\dlog{1}{1}$.

\begin{proposition} \label{tr}
Let $A=\dlog{\wek s}{\wek p}$ be a matrix coded by a double 0-1 sequence of depth $n$.

\begin{enumerate}
\item
 The transpose $\ct A=(-1)^{{{c_{11}}}} A$.
\item
$\tr \dlog{\wek 0}{\wek 0}=2^n$ and $\tr \dlog{\wek s}{\wek p}=0$ iff $\dlog{\wek s}{\wek p}\neq \dlog{\wek 0}{\wek 0}$.
\item
If $B$ is coded by a distinct double sequence, then $\tr (A^TB)=0$.
\end{enumerate}
\end{proposition}
\Proof
The first two statements follow by simple mathematical induction with respect to depth, and the third is their consequence.
\QED
\vspace{3pt}

The eigen-structure of doubled matrices also follows a simple doubling algorithm.
\begin{proposition}
For an $n\times n$ matrix $e$ with an eigenvalue $\lambda$ corresponding to an eigenvector $\wek x$, the doubled matrix $\dlog s p (e)$ entails the following pattern of new eigenvectors and eigenvalues:
\def\bpmat{\begin{pmatrix}}
\def\epmat{\end{pmatrix}}
\[
\begin{array}{lcrccrcc}
\vspandex
\displaystyle\dlog 0 0(e)=
\bpmat\hh e\hh&\hh0\hh\\0&e\epmat&\qquad\mapsto\qquad&\lambda&\sim
          &\bpmat \hhh\wek x\hhh\\\hhh\wek 0\hhh\epmat,&\lambda&\sim
          &\bpmat \hhh\wek 0\hhh\\\wek x\epmat;\\
\vspandex
\displaystyle\dlog 1 0(e)=
\bpmat\hh e\hh &\hh 0\hh\\0&-e\epmat&\qquad\mapsto\qquad& \lambda&\sim
          &\bpmat \hhh\wek x\hhh\\\hhh\wek 0\hhh\epmat,&-\lambda&\sim
          &\bpmat \hhh\wek 0\hhh\\\hhh\wek x\hhh\epmat;\\
\vspandex
\displaystyle\dlog 0 1(e)=
\bpmat\hh 0\hh &\hh e\hh\\e&0\epmat&\qquad\mapsto\qquad& \lambda&\sim
          &\bpmat \hhh\wek x\hhh\\\hhh\wek x\hhh\epmat,&-\lambda&\sim
          &\bpmat \hhh\wek x\hhh\\-\wek x\epmat;\\
\vspandex
\displaystyle\dlog 1 1(e)=
\bpmat\hhh 0\hh &\hh e\hh\\-e&0\epmat&\qquad\mapsto\qquad& \imath\lambda&\sim
          &\bpmat\hhh\wek x\hhh\\-\imath\wek x\epmat,&-\imath\lambda&\sim
          &\bpmat \hhh\wek x\hhh\\\imath\wek x\epmat.
\end{array}
\]
Consider $\dlog{\wek s}{\wek p}\neq \dlog{\wek 0}{\wek 0}$. In particular, its eigenvalues are either $\pm 1$ when ${{c_{11}}}$ is even, or $\pm \imath$ otherwise. Also, the algorithm preserves orthogonal eigenvectors.
\end{proposition}
\Proof
The doubling of eigenvalues is seen through the characteristic polynomials, which  are, consecutively:
 \[
 |e-\lambda I|^2,\quad
 (-1)^n|e-\lambda I|\cdot|e+\lambda I|,\quad
 (-1)^n|e-\lambda I|\cdot|e+\lambda I|,\quad
 |e+\imath \lambda I|\cdot|e-\imath \lambda I|.
 \]
The eigenvectors as well as the corollaries follow directly.
\QED\vspace{3pt}

Now, we will focus on the embedding aspects \refrm{Qact} of operators $Q^{\wek s}_{\wek p}:\R^{k\times k}\to \R^{2^nk\times 2^nk}$, where $n$ is the depth of the double sequence.
\begin{proposition}
Denote the entries of the matrix $Q=Q^{\wek s}_{\wek p}$ by $q_{jk}$. Let a ${2^nk\times 2^nk}$ matrix $A$ be partitioned into ${2^{2k}}$ blocks $A_{jk}$ of size $n\times n$ The adjoint operators $P=Q^*:\R^{2^nk\times 2^nk}\to \R^{k\times k}$ act as follows:
\[
PA=\frac{1}{2^n}\sum_{j,k} q_{jk} A_{jk}.
\]
In particular, the principal quadrupling $\wek s=\wek p=\wek 0$ yields the trace-like adjoint
\[
PA=\frac{1}{2^n}\sum_{j} A_{jj}.
\]
\end{proposition}
\Proof
We check the formulas directly first for $n=1$, and then proceed by induction.\QED
\subsection{Spanned algebras}

While the group $\ca D$ does not have a natural tangible extension to an algebra yet its isomorphic copy $\ca C$ extends to the operator algebra  $\cC$ (defined below \refrm{action}). Alternatively, its isomorphic image $\ca M$ extends to the matrix algebra $\cM$ (defined below \refrm{action1}).
\vv

A vector space of finitary real matrices can be turned into an inner product space by putting
\[
\la K,M\ra\df \tr \left(\ct K M\right),
\]
entailing the Hilbert-Schmidt (or Frobenius) norm.
However, the matrices represented by  double 0-1 finitary sequences will have norms escaping to infinity along with the increasing depth. Therefore, at least initially, we need to curtail the depth, confining ourselves to finite dimensional spaces. So, consider the $2^{2n}$-dimensional real Hilbert space of matrices coded by distinct double sequences $\dlog {\wek s}{\wek p}$ with the maximum depth $n$, with the orthonormal basis (cf. Proposition \ref{tr}.2)
\[
\frac{1}{2^{n/2}}\dlog{\wek s}{\wek p}.
\]
In particular, a basic sequence of length $k\le n$, with $n$ being the maximum of the depths of the sequence's elements, spans  the algebra, equipped with the above inner product.

\begin{example}[Clifford algebras]\rm
For $k\le n$, consider basic $e_1,\dots, e_k$ selected from $\dlog {\wek s}{\wek p}\neq \dlog{\wek 0}{\wek 0}$ of depth $n$. Assume that the generator is {\em pure}, i.e., all members have the same sign, either all positive ($+$) or all negative ($-$). Then $E=\R^{2^k}$ can be embedded into the algebraic span $\cA$ of the sequence through the formula
\[
j(\wek x)=\sum_{\wek d} x_{\wek d}\, \wek e^{\wek d},\quad \wek d\in\set{0,1}^n, \quad \wek x=(x_{\wek d}).
\]
$E$ is naturally equipped with the quadratic form
\[
q(\wek x)=\frac{1}{2^n} \sum_{\wek d} |x_{\wek d}|^2.
\]
Clearly, by Proposition \ref{tr},
\vv
\vspandexsmall \hspace{20pt} (i) $1=\dlog{\wek 0}{\wek 0} \notin j(E)$, since $\dlog{\wek 0}{\wek 0}$ is orthogonal to all powers $\wek e^{\wek d}$.

\vspandexsmall\hspace{20pt} (ii) $(j(\wek x))^2=q(\wek x) $ (under ($+$)) or $(j(\wek x))^2=-q(\wek x) $ (under ($-$)) for every $\wek x$,

\vspandexsmall \hspace{20pt} (iii)  $j(E)$ spans $\cA$.

\vv
In other words, $(\cA,j)$ becomes a Clifford algebra for $(E,q)$ (cf. \cite[Sect.\ 5.1]{Gar} or similar texts).
\end{example}

\section{A rigged Hilbert space}\label{RHS}
\subsection{A directed injective system of Hilbert spaces}
In doubling (or quadrupling) we recognize a directed system of $C^*$-algebras that also happen to be Hilbert spaces of finite dyadic dimension. The essential general framework can be found, e.g., in \cite[Prop. 11.4.1 and Exerc. 11.5.26]{KadRin} and further details in the context of Hilbert spaces -  in \cite{BelTra}.\vv

\def\uu#1#2{{_#2U^#1}}
Let a double sequence $\dlog{\wek s}{\wek p}$ have depth $n$. As a matrix, it is a member of the real Hilbert space and algebra $\ca H_n=\R^{2^n\times 2^n}$, equipped with the inner product
\[
\la A,B\ra_n=\frac{1}{2^n}\,\tr
(\ct AB)
\]
In the isomorphic context of Rademacher chaos operators, the domain of operators $\wek R^{\wek s}\wek S^{\wek p}$ is restricted to $L^2(\ca D_n)$ (cf. the vicinity of \refrm{Doob}), equipped with the normalized Hilbert-Schmidt norm.\vv

By convention, $\ca H_0=\R$. Then
$2^{2n}$ of matrices coded by distinct double sequences form an orthonormal basis of $\ca H_n$ (cf.\ Proposition \ref{tr}). The principal quadrupling operator
 \[
 Q=Q_0^0,\quad Q^2=Q_{00}^{00}, \quad Q^3=Q_{000}^{000},\quad ...
 \]
 (stretching somewhat the power notation for the sake of transparency) entails the consistent family of isometric injections between the members of the linearly ordered sequence of Hilbert spaces:
\[
\ca H_n\stackrel{\uu n m}{\longrightarrow} \ca H_m:\qquad \uu n m A\df Q^{m-n} A,\quad m\ge n, \quad A\in\ca H_n.
\]
By convention, $\uu n n $ is the identity mapping and by consistency we understand the property
\[
\uu n k = \uu m k \uu n m, \quad k\ge m\ge n.
\]
\begin{remark}\rm
We will not elaborate on the injections. However, if we did, then it would be convenient to use the left subscript for  members of our spaces as in the formula ${_m x}= \uu  n m{_n x},\,m\ge n$.
\end{remark}
\vv
Denote $\BH=\Prod_{n\ge 0} \ca H_n $. The $C^*$-algebra $\ell^\infty(\BH)$ of bounded sequences, equipped with the supremum norm, and its closed subalgebra $c_0(\BH)$ of null sequences entail the quotient $C^*$-algebra
\[
\widetilde{\cS}=\ell^\infty(\BH)/c_0(\BH)=\set{\widetilde{\wek A}=\wek A+c_0(\BH):\wek A\in \ell^\infty(\BH) }.
\]
The isometric injection
\be\label{QA}
A_k\stackrel{\theta_k}{\longrightarrow}(\underbrace{0,...,0}_{k-1},A_k,Q A_k,Q^2 A_k,\dots),
\ee
extends to the isometric injection $\widetilde{\theta}_k\wek A=\theta\wek A+c_0(\BH)$ from $\ca H_k$ into $\bigcup_k \widetilde{\theta}_k(\ca H_k)$, the latter being endowed with the inner product
\be\label{AB}
\la \widetilde{\wek A},\widetilde{\wek B}\ra =\lim_k \la A_k,B_k\ra_k,
\ee
(the sequence on the right is constant eventually). The completion under the derived norm becomes a Hilbert space into which all $\ca H_n$ embed isometrically.
\subsection{An alternative}
The construction of the Hilbert space, presented in \cite{BelTra}, departs to some extent from the universal framework of \cite[Prop.\ 11.4.1]{KadRin}. Using a simpler approach the authors deal with a system of Hilbert spaces  directed by contractions rather than isometries. \vv

First, their starting vector space $\cS\subset\BH $ consists of sequences $\wek A=(A_n)$ such that for some $k$ and for all $j\ge 0$, $A_{k+j}=Q^jA_k$. Next,
the stronger equivalence is considered: two sequences $\wek A$ and $\wek B$ are deemed equivalent if their elements are equal eventually (as opposed to differ by a null sequence with respect to the supremum norm). We will keep the same notation though: $\wek A\sim\wek B$. More precisely,
the vector subspace
$\cS_0=\set{\wek A: A_k =0\mbox{ eventually}}$ entails the quotient space $\widetilde{\cS}=\cS /\cS_0$. Under the inner product \refrm{AB} the injections $\widetilde{\theta_n}:\ca H_n\to \widetilde{\cS}$ are isometries, and their ascending ranges yield the inductive limit. Under the notation of  \cite{BelTra},
\[
\cD^\times\df \lim_{\longrightarrow} \widetilde{\theta}_k(\ca H_k),
\]
which is a vector subspace of $\widetilde{\cS}$ equipped with the strongest locally convex topology ensuring the continuity of all inclusions.
A typical equivalence class consists of sequences of the following form, for some $k\ge 0$ and $A_j\in \ca H_j,\, j=1,\dots,k$:
\[
(A_0,A_1,\dots, A_{k-1},A_k,Q A_k,Q^2 A_k,\dots).
\]
Let us use the acronym NQ for a matrix that is not quadrupled, i.e., a matrix $A\neq QB$, for any matrix $B$ of half-size of $A$. For each sequence $\wek A\in\cS$ or $\widetilde{\wek A}=\wek A+\cS_0\in\widetilde{\cS}$ we can define its {\em quadrupling index} and extract the first (with respect to the natural order) NQ element:
\[
d(\wek A)=\min\set {k: A_j=Q^{j-k} A_k,\,j\ge k}, \quad\mbox{yielding} \quad A_{d(\wek A)}.
\]
Naturally, $d(\wek A+\cS_0)=\Min\set{d(\wek A)}$. Although the notion of ``the first NQ element'' is well defined yet its usage is rather limited, since for given $\wek A\in\cS$
\[
\forall\,\epsilon>0\quad \sup\set{d(\wek B): ||\wek A-\wek B||<\epsilon}=\infty.
\]
Indeed, let $k=d(\wek A)$ and $A_k$ be the first NQ element. Then define $\wek B^{(j)}$ through its first NQ element
$B_{k+j}=Q^jA_k+R_{k+j}$, where $||R_{k+j}||<\epsilon$. In particular, it is possible that $\wek B^{(j)}\to \wek A$ yet $d(\wek B^{(j)})\to\infty$. Nevertheless, this shows that $\widetilde{\cS}$ is not complete since $\wek A+c_0(\BH)$ contains cluster points of $\wek A+\cS_0$.

\vv
On the other hand the notion of ``NQ'' sheds some light on the role of the adjoint operator $P=Q^*$:
\[
\ddmat {A_{11}}{A_{12}}{A_{21}}{A_{22}} \quad\stackrel{P}{\longrightarrow} \quad \frac{1}{2}\left(A_{11}+A_{22}\right).
\]
\begin{proposition} Fix $k\ge 2$ and let a NQ $A\in \ca H_k$ generate the sequence $\theta_k(A)\in \cS$. Let $j<k$. Then the infimum in the distance formula
\[
{\rm dist} \left(\theta_j(\ca H_j),\theta_k(A)\right)=\inf\set{||Q^{k-j}B-A||:B\in \ca H_j}
\]
is attained at $B=P^{k-j}A$.
\end{proposition}
\Proof
Let $A$ be partitioned into $2^{2(k-j)}$ blocks $A_{lm}$ of size $2^j\times 2^j$. Then
\[
||Q^{k-j}B-A||^2=\frac{1}{2^{k-j}} \left(\sum_l ||A_{ll}-B||^2+\sum_{l\neq m} ||A_{lm}||^2  \right)
\]
 We are dealing with the classical regression (or projection) problem in a Hilbert space. That is, for a random variable $X$ with values in a Hilbert space $H$, we have
\[
\inf_{x\in H} \E||X-x||^2=\E ||X-\E X||^2,
\]
and in our case we have just a random variable $X$ uniformly distributed on the finite set $\set{A_{ll}}$. Clearly, $\E X=P^{k-j} A$.
\QED
\vspace{3pt}

The observation leads to the candidate for the dual space of $\cD^{\times}$. Under the notation of \cite{BelTra}
\[
\cD\df\set{(A_k): A_j=P^{k-j} A_k, j\le k},
\]
where the duality $(\cD,\cD^\times)$ is defined by the bilinear form \refrm{AB}. Again, the scalar sequence is constant eventually. $\cD$ is equipped with the product topology $\tau_p$, i.e., the weakest topology that ensures the continuity of the surjective projections $\pi_k:\cD\to \cH_k$. In fact, they are contractions. Then, indeed, the dual of $(\cD,\tau_p)$ is $\cD^\times$ (\cite[Theorem 2.2 and Lemma 2.1 (iii)]{BelTra}). \vv
\begin{remark}\rm
In \cite{BelTra} the adjoint operators $(_mU^n)^*$ were assumed to be one-to-one, which does not occur in our case since the adjoint operators are projections. The injectivity  could be ensured, e.g.,  by dense ranges of embeddings ${_mU^n}$, motivated by specific applications in the paper (C.\ Trapani, personal communication, July 31, 2016). Fortunately,  that does not affect the derived important properties of the triple $(\cD,(\ca H_n),\cD^\times)$ such as duality. The triple is named {\em the rigged Hilbert space}.
\end{remark}
\section{Summary}
From the algebraic point of view, a Rademacher system is generated by a commutative signed positive sequence $\wek r=(r_n)$ in a Hilbert space such that the powers $\wek r^{\wek p}, \,\wek p\in \BD,$ form an orthonormal sequence, yielding the real Hilbert space and algebra, followed by analysis  of counterparts of differential operators, called here ``chaos operators''. The study of these operators is motivated by applications in quantum probability. As a part of the analysis we have derived forms of the Leibnitz Formula and the Chain Rule.\vv

 The commutative patterns of these operators and their compositions are fairly complex with exception of two operators whose products and powers algebraically span all other operators. Their actions are recorded as operations on double binary sequences and the obtained multiplicative system can be represented by matrices of dyadic size, i.e., $2^n\times 2^n$, for some $n$. \vv

 In virtue of Theorem \ref{embed} and in view of the construction of a rigged Hilbert space in the latter section, we can create a variety of structures analogous to Rademacher systems.  They share the Hilbertian and algebraic structure with an orthonormal basis of powers. At the same time, their signs and commutative properties can be chosen at will as in Example \ref{examples:embed}.\vv

 Since the constructed (or perhaps only recorded) structures extend the notions of Pauli spin matrices and quaternions, they may find potential applications in quantum theories. In addition, due to the binary and algorithmic character of the constructions, computer simulations may shed new light on chaotic behavior.

\end{document}